\documentclass[11pt]{article}
\usepackage{amscd}
\usepackage{amsfonts}
\usepackage{amsmath}
\usepackage{amssymb}
\usepackage{amsthm}
\usepackage{bbm}
\usepackage{CJK}
\usepackage{fancyhdr}
\usepackage{graphicx}
\usepackage{hyperref}
\usepackage{indentfirst}
\usepackage{latexsym}
\usepackage{mathrsfs}
\usepackage{xypic}

\usepackage[top=1in,bottom=1in,left=1.25in,right=1.25in]{geometry}
\def\<{\langle}
\def\>{\rangle}
\def\a{\alpha}
\def\b{\beta}

\def\c{\cdot}

\def\o{\otimes}

\date{}
\begin{document}
\renewcommand{\baselinestretch}{1.2}
\renewcommand{\arraystretch}{1.0}
\title{\bf  On BiHom-analogue of generalized Lie  algebras}
\author{{\bf Shuangjian Guo$^{1}$, Xiaohui Zhang$^{2}$,  Shengxiang Wang$^{3}$\footnote
        { Corresponding author(Shengxiang Wang):~~wangshengxiang@chzu.edu.cn} }\\
{\small 1. School of Mathematics and Statistics, Guizhou University of Finance and Economics} \\
{\small  Guiyang  550025, P. R. of China} \\
{\small 2.  School of Mathematical Sciences, Qufu Normal University}\\
{\small Qufu  273165, P. R. of China}\\
{\small 3.~ School of Mathematics and Finance, Chuzhou University}\\
 {\small   Chuzhou 239000,  P. R. of China}}
 \maketitle

\begin{center}
\begin{minipage}{13.cm}

{\bf \begin{center} ABSTRACT \end{center}}
 In this paper, we introduce the definition of generalized BiHom-Lie  algebras and generalized BiHom-Lie admissible algebras
in the category ${}_H{\mathcal M}$ of left modules for any quasitriangular Hopf algebra $(H, R) $.
 Also, we  describe the BiHom-Lie ideal structures of the BiHom-associative algebras.\\


 \end{minipage}
 \end{center}
 \normalsize\vskip1cm

\section*{Introduction}
\def\theequation{0. \arabic{equation}}
\setcounter{equation} {0}

Hom-algebras were firstly studied by Hartwig, Larsson and Silvestrov \cite{Hartwig},
where they introduced the structure of Hom-Lie algebras in the context of the deformations of
Witt and Virasoro algebras.
Later, motivated by the new examples arising as applications of the general quasi-deformation construction
and the desire to be able to treat within the same framework of the super and color Lie algebras,
Larsson and Silvestrov \cite{Larsson} extended the notion of Hom-Lie algebras to quasi-Hom Lie algebras
and quasi-Lie algebras, such that the classical definition of this algebraic
structure is "deformed" by means of this endomorphism. The theory of Hom-type algebraic
structures has seen an enormous growth in recent years.

An elementary but important property of Lie algebras is that each associative algebra
gives rise to a Lie algebra via the commutator bracket.
Makhlouf and Silvestrov (\cite{Makhlouf2008}, \cite{Makhlouf2009}, \cite{Makhlouf2010}) generalized the associativity to twisted associativity
and naturally proposed the notion of Hom-associative algebras.
Furthermore they obtained that a Hom-associative algebra gives rise to a Hom-Lie
algebra via the commutator bracket. Caenepeel and Goyvaerts \cite{Caenepeel2011} studied the Hom-Hopf algebras from a categorical view point. In nowadays mathematics, much of the research on certain algebraic object is to
study its representation theory. The representation theory of an algebraic object reveals some of its profound structures hidden underneath. A good example is
that the structure of a complex semi-simple Lie algebra is much revealed via its
representation theory. Sheng \cite{Sheng2012} studied the representation theory of Hom-Lie algebras.  Later on, Wang \cite{Wang2014}  studied the representation theory of Hom-Lie algebras in Yetter-Drinfeld categories.

A BiHom-algebra is an algebra in such a way that the identities defining the structure are
twisted by two homomorphisms $\a,\b$. This class of algebras was introduced from a categorical
approach in \cite{Graziani} as an extension of the class of Hom-algebras. More applications of BiHom-Lie
algebras, BiHom-algebras, BiHom-Lie superalgebras and BiHom-Lie colour algebras
can be found in (\cite{Abdaoui2017}, \cite{Cheng2016}, \cite{lijuan}, \cite{Wang2016}, \cite{zhang18}). What forms in the BiHom-case are
Wang's results in \cite{Wang2014}? In this paper, we give a positive answer to the question.

This article is organized as follows.  In Section 2, we introduce the notion of  BiHom-Lie  algebras
in ${}_H{\mathcal M}$ and show that a BiHom-associative algebra in ${}_H{\mathcal M}$ gives rise to a
BiHom-Lie algebra in ${}_H{\mathcal M}$ by the natural bracket product (see Theorem 2.3).
In Section 3, we determine the the BiHom-Lie ideal structures of the BiHom-associative algebras
(see Theorem 3.9).

\section{Preliminaries}
\def\theequation{\arabic{section}.\arabic{equation}}
\setcounter{equation} {0}

In this section we recall some basic definitions and results related to our paper.
Throughout the paper, all algebraic systems are supposed to be over a field ${k}$.
The reader is referred to \cite{Graziani}
as general references about BiHom-associative algebras.
If $C$ is a coalgebra, we use the Sweedler-type notation \cite{Sweedler} for the comultiplication: $\Delta(c)=c_{1}\o c_{2}$,
for all $c\in C.$

Recall that if  $H$ is a  bialgebra and $M$ is a left $H$-comodule with coaction
\begin{eqnarray*}
\rho: H\o M\rightarrow M, ~~~~h\o m\mapsto  h\c m,~~~\forall h\in H,  m\in M,
\end{eqnarray*}
the coassociativity of the coaction means $(\Delta\o id )\circ \rho=(id \o \rho)\circ \rho$.

In this paper, we consider objects in the category of left $H$-modules ${}_H{\mathcal M}$. In
particular, a left $H$-module algebra $A$ is an algebra in this category; this means
that multiplication in $A$ is an $H$-module map:
\begin{eqnarray}
h\c (ab)=(h_1\c a)(h_2\c b),~~~\forall a,b\in A.
\end{eqnarray}

\noindent{\bf 1.1. Quasitriangular structure} A pair $(H, R)$ is called a  quasitriangular Hopf algebra if $H$ is a Hopf algebra and $R=\sum R^{1}\o R^2\in H\o H$   of invertible elements (called a universal $R$-matrix), satisfying the following conditions:
\begin{eqnarray*}
&&(QT1) (\Delta \o H)(R) = R_{13} R_{23},\\
&&(QT2) (H \o \Delta)(R) = R_{13} R_{12},\\
&&(QT3) R\Delta(h) = \Delta^{cop}(h)R, 
\end{eqnarray*}
for any $h\in H$. If, in addition, $R$ is symmetric, that is
\begin{eqnarray*}
(e)&& R^{-1}=\sum R^{2}\o R^1,
\end{eqnarray*}
then $(H,R)$ is called triangular.

In the category ${}_H{\mathcal M}$, the braiding $\tau: M\o N\rightarrow N\o M$
is given by
 $$
 \tau(m\o n) =\sum R^2\c n \o  R^1\c m,
 $$
for all $m \in M \in {}_H{\mathcal M}$ and $n \in N \in {}_H{\mathcal M}$.

Let $A$ be an algebra in ${}^H{\mathcal M}$,
$A$ is called $H$-commutative if
 \begin{eqnarray}
\sum(R^2\c b) \o  (R^1\c a)=ab,
\end{eqnarray}
for all $a,b\in A.$

\noindent{\bf 1.2. BiHom-associative algebra} Recall from \cite{Graziani}
that a BiHom-associative algebra  is a 4-tuple $(A,m,\alpha,\b)$
consisting of a linear space $A$, a bilinear map $m:A\otimes A\rightarrow A$
and  homomorphisms $\alpha, \b: A\rightarrow A$ such that, for all $a,b,c\in A$,
 \begin{eqnarray}
 \a\circ\b=\b\circ\a,~~
\alpha(a)(bc)=(ab)\b(c).
\end{eqnarray}

In particular, if $\alpha(ab)=\alpha(a)\alpha(b)$ and $\b(ab)=\b(a)\b(b)$, we call  $A$ a multiplicative BiHom-associative algebra.
If there exists an element $1_{A}\in A$ such that
$1_{A}a=\b(a)$ and $a1_{A}=\alpha(a)$ for all $a\in A$,
we call $A$  a  unital BiHom-associative algebra.

\noindent{\bf 1.3. BiHom-Lie algebra } Recall from \cite{Graziani}
 that  a BiHom-Lie algebra is a 4-tuple  $(L, [,],\alpha,\b)$ consisting of a linear space $L$,
a bilinear map $[,]:L\o L\rightarrow L$
and  homomorphisms $\alpha, \b: L\rightarrow L$ satisfying:
\begin{eqnarray*}
&&  \a\circ\b=\b\circ\a,\\
&&[\b(l),\a(l')]=-[\b(l'),\a(l)],~(Skew\mbox{-}symmetry),\\
&&\circlearrowleft_{l,l',l''}[\b^2(l),[\b(l'),\a(l'')]]=0,~(BiHom\mbox{-}Jacobi~identity),
\end{eqnarray*}
 for any  $l,l',l''\in L$, where $\circlearrowleft$ denotes the summation over the cyclic permutation
on  $l,l',l''.$

\section{Generalized BiHom-Lie  algebras}
\def\theequation{\arabic{section}.\arabic{equation}}
\setcounter{equation} {0}

In this section we introduce the concept of BiHom-Lie algebras
in the category ${}_H{\mathcal M}$ and
provide a construction of BiHom-Lie  algebras in ${}_H{\mathcal M}$ through BiHom-associative algebras in ${}_H{\mathcal M}$.
\medskip

\noindent{\bf Definition 2.1.}
Let $(H, R)$ be a quasitriangular Hopf algebra.
A BiHom-Lie algebra in the category ${}_H{\mathcal M}$
(called  a generalized  BiHom-Lie algebra)
 is a 4-tuple $(L,[,],\alpha,\b)$,
where $L$ is an object in ${}_H{\mathcal M}$,
$\alpha, \b: L\rightarrow L $ are  homomorphisms in ${}_H{\mathcal M}$
and $[,]:L\o L\rightarrow L$ is a morphism in ${}_H{\mathcal M}$ satisfying:
\begin{eqnarray}
 &&\a\circ\b=\b\circ\a,\\
 &&\a([l,l'])=[\a(l),\a(l')],~~~\b([l,l'])=[\b(l),\b(l')],\\
 &&[\b(l),\a(l')]=- [R^2\c \b(l'),R^1\c \a(l)],\\
 &&\{l\o l'\o l''\}+\{(\tau\o 1)(1\o\tau)(l\o l'\o l'')\}+\{(1\o\tau)(\tau\o 1)(l\o l'\o l'')\}=0,~~~~
\end{eqnarray}
 for any  $l,l',l''\in L$, where $\{l\o l'\o l''\}$ denotes
$[\b^2(l),[\b(l'),\a(l'')]]$ and $\tau$ the braiding for $L$.
\medskip

\noindent{\bf Example 2.2.}
 If $\b=\a$ in Definition 2.1, then the generalized  BiHom-Lie algebra $L$
is just the generalized Hom-Lie algebra in Wang \cite{Wang2014}.
\medskip

\noindent{\bf Theorem 2.3.}
Let $(H, R)$ be a triangular Hopf algebra and $(A,\alpha,\b)$ a BiHom-associative algebra in ${}_H{\mathcal M}$
(called  a generalized  BiHom-associative algebra) with two bijective homomorphisms $\a$ and $\b$. Then the 4-tuple $(A,[,],\alpha,\b)$ is a generalized BiHom-Lie algebra,
where the bracket product $[,]:A\o A\rightarrow A$ is defined by
\begin{eqnarray}
[a,b]=ab- (R^2\c \a^{-1}\b(b))(R^1\c \a\b^{-1}(a)),
\end{eqnarray}
for all $a,b\in A$.

{\it Proof } Similar to \cite{Abdaoui2017}.
\medskip

\noindent{\bf Example 2.4.}
Let $\{x_{1},x_{2}\}$ be a basis of a 2-dimensional linear space $A$. The
following multiplication $m$ and linear maps $\alpha,\b$ on $A$ define a BiHom-associative algebra (\cite{Graziani}) :
\begin{eqnarray*}
&&m(x_{1},x_{1})=x_{1},m(x_{1},x_{2})=bx_{2},\\
&&m(x_{2},x_{1})=-x_{2},m(x_{2},x_{2})=0,\\
&&\alpha(x_{1})=x_{1},\alpha(x_{2})=-x_{2},\\
&&\b(x_{1})=x_{1},\b(x_{2})=bx_2,
\end{eqnarray*}
where $b$ is a  parameter in $k$.

Let $G$ be the cyclic group of order $2$ generated by $g$.
The group algebra $H=kG$ is a Hopf algebra in the usual way.
 We take $R =\frac{1}{2} (e\o e+e\o g +g \o e-g \o g)$,
where $e$ is the unit of the group $G$. It is easy to check that $H$ is a triangular Hopf algebra.
Define the left-$H$-comodule structure of $A$ by
\begin{eqnarray*}
&&e\c x_{i})=x_{i},~~~g\c x_{1}=x_{1}, g\c x_{2}=-x_{2}
\end{eqnarray*}
It is not hard to check that $A$ is a generalized BiHom-associative algebra.

The braiding $\tau$ is given by $\tau(x_{2}\o x_{2})=-x_{2}\o x_{2}$
and $\tau(x_{i}\o x_{j})=x_{j}\o x_{i}$ for other cases.
Then according to Theorem 2.3,
we obtain a generalized BiHom-Lie algebra with the bracket product
$[,]$ satisfying the following non-vanishing relation
\begin{eqnarray*}
[x_{1},x_{2}]=2bx_2,~~~   [x_{2},x_{1}]=-x_{2}+bx_{1}.
\end{eqnarray*}

\noindent{\bf Example 2.5.}
Let $A$ be the three dimensional Heisenberg Lie algebra,
which consists of the strictly upper-triangular complex $3\times3$ matrices.
It has a standard linear basis
 \[
 \left(
   \begin{array}{ccc}
     0 & 1 & 0 \\
     0 & 0 & 0 \\
     0 & 0 & 0 \\
   \end{array}
 \right),
 \left(
   \begin{array}{ccc}
     0 & 0 & 0 \\
     0 & 0 & 1 \\
     0 & 0 & 0 \\
   \end{array}
 \right),
 \left(
   \begin{array}{ccc}
     0 & 0 & 1 \\
     0 & 0 & 0 \\
     0 & 0 & 0 \\
   \end{array}
 \right).
 \]

Let $G$ be the cyclic group of order $2$ generated by $g$.
The group algebra $H=kG$ is a Hopf algebra in the usual way, we take 
$R =\frac{1}{2} (e\o e+e\o g +g \o e-g \o g)$, it is easy to check that $H$ is a triangular Hopf algebra.  Define the left-$H$-module structure of $A$  by
\begin{eqnarray*}
&&e\c e_i=e_i,  g\c x_{1}=-x_1, g\c x_{2}=x_{2},g\c x_{3}=x_{3},
\end{eqnarray*}
where $e$ is the unit of the group $G$.
It is not hard to check that $A$ is an object in ${}_H{\mathcal M}$.

The braiding $\tau$ is given by
\begin{eqnarray*}
&&\tau(x_{1}\o x_{1})=-x_{1}\o x_{1},\tau(x_{1}\o x_{2})=-x_{2}\o x_{1},\tau(x_{1}\o x_{3})=x_{3}\o x_{1},\\
&&\tau(x_{2}\o x_{1})=-x_{1}\o x_{2},\tau(x_{2}\o x_{2})=-x_{2}\o x_{2},\tau(x_{2}\o x_{3})=x_{3}\o x_{2},\\
&&\tau(x_{3}\o x_{1})=x_{1}\o x_{3},~~\tau(x_{3}\o x_{2})=x_{2}\o x_{3},~~\tau(x_{3}\o x_{3})=x_{3}\o x_{3}.
\end{eqnarray*}

Define a bracket product $[,]$ on $A$ by
\begin{eqnarray*}
&&[x_{i},x_{3}]=[x_{3},x_{i}]=[x_{1},x_{1}]=[x_{2},x_{2}]=0,\\
&&[x_{1},x_{2}]=[x_{2},x_{1}]=x_{3},~i=1,2,3.
\end{eqnarray*}
One may verify that $(A,[,])$ is a generalized Lie algebra.

Let $\lambda_{1},\lambda_{2}$ be two nonzero scalars in $k$.
Consider the maps
$\alpha,\b:A\rightarrow A$
defined on the basis element by
\begin{eqnarray*}
\alpha(x_{1})=\lambda_{1}x_{1},
\alpha(x_{2})=\lambda_{2}x_{2},
\alpha(x_{3})=\lambda_{1}\lambda_{2}x_{3},\\
\b(x_{1})=\lambda'_{1}x_{1},
\b(x_{2})=\lambda'_{2}x_{2},
\b(x_{3})=\lambda'_{1}\lambda'_{2}x_{3}.
\end{eqnarray*}
It is straightforward to check that $\alpha,\b$ defines
a Lie algebra endomorphism in ${}_H{\mathcal M}$.
We  obtain a generalized BiHom-Lie algebra $(A,[,]',\alpha,\b)$, whose bracket product
satisfies the following non-vanishing relation
\begin{eqnarray*}
[x_{1},x_{2}]'=\lambda_{1}\lambda'_{2}x_{3},~~~[x_{2},x_{1}]'=\lambda'_{1}\lambda_{2}x_{3}.
\end{eqnarray*}

\section{On the Generalized BiHom-Lie ideals
structures of generalized BiHom-associative algebras}
\def\theequation{\arabic{section}.\arabic{equation}}
\setcounter{equation} {0}

In this section,
we consider some $H$-analogous of classical concepts of ring theories
and  Lie theories as follows.
\medskip

Let $(A,\a,\b )$ be a generalized BiHom-associative algebra.
An $H$-BiHom-ideal $U$ of $A$ is an ideal
such that $\alpha(U)\subseteq U, \b(U)\subseteq U$ and $(AU)A=A(UA)\subseteq U.$
\medskip

Let $(L, \a, \b)$ be a generalized BiHom-Lie algebra.
An  $H$-BiHom-Lie ideal $U$ of $L$ is a Lie ideal
such that $\alpha(U)\subseteq U, \b(U)\subseteq U $ and $[U,L]\subseteq U.$
Define the  center of $L$ to be $Z_{H}(L)=\{l\in L|[l,L]_{H}=0\}.$
\medskip

$(L, \a, \b)$ is called  prime if the product of any two non-zero $H$-BiHom-ideals of $L$ is non-zero.
It is called  semiprime if it has no  non-zero nilpotent $H$-BiHom-ideals,
and is called simple if it has no nontrivial $H$-BiHom-ideals.
\medskip

\noindent{\bf Lemma 3.1.}
Let $(A,\alpha, \b)$ be a generalized  BiHom-associative algebra. Then

(1) $[\alpha\b(a),bc]=[\b(a),b]\b(c)+ R^2\c \b(b)[ R^1\c\a(a),c],$

(2) $[ab,\alpha\b(c)]=\a(a)[b,\a(c)]+ [a,  R^2\c\b(c_0)] R^1\c\a(b_{0}),$

\noindent for all $a,b,c\in L.$

{\it Proof }  We only prove (1), and similarly for (2). For any $a,b,c\in L,$
\begin{eqnarray*}
&&[\b(a),b]\b(c)+ R^2\c\b(b)[ R^1\c\a(a),c]\\
&=&(\b(a)b)\b(c)-((R^2\c \a^{-1}\b(b))(R^1\c\a\b^{-1}\b(a)))\b(c)+(R^2\c \b(b))(R^1\c \a(a)c)\\
&&-(R^2\c \b(b))((r^2\c\a^{-1}\b(c))(r^1R^1\c \a^{2}\b^{-1}(a)))\\
&=&\alpha\b(a)(bc)-(R^2_1\c \b(b))((R^2_2\c\a^{-1}\b(c))(R^1\c \a^{2}\b^{-1}(a)))\\
&=&\alpha\b(a)(bc)-R^2\c\a^{-1}\b(bc)(R^1\c \alpha\b^{-1}(\a\b(a))\\
&=&[\alpha\b(a),bc].
\end{eqnarray*}

\noindent{\bf Lemma 3.2.}
Assume that $(L, \a,\b)$ is a semiprime unital generalized BiHom-associative algebra,
and let $U$ be both an $H$-BiHom-Lie ideal and a BiHom-associative subalgebra of $L$.
If $[U,U]\neq 0,$
then there exists a non-zero $H$-BiHom-ideal of $L$.
\medskip

{\it Proof }  Since $[U,U]\neq 0,$ there exists $x,y\in U$ such that $[x,y]\neq 0$,
and $[lx,\alpha(y)]\in U$ for all $l\in L$.
By Lemma 3.1(2),
$\a(l)[x,\a(y)]=[lx,\alpha\b(y)]-\sum[l, R^2\c \b(y)]R^1\c \a(x)$,
thus $[l,  R^2\c \b(y)]\in U$.
Because $U$ is an BiHom-associative subalgebra of $L$,
then $[l,   R^2\c \b(y)]R^1\c \a(x)\in U,$
therefore $\a(l)[x,\a(y)]\in U$ for all $l,m\in L.$
It follows that
$I=(L[x,y])L=L([x,y]L)\in U.$
In fact, take $I=\{\sum_{i}(a_{i}\alpha^{n_{i}}\b^{m_i}([x,y]))b_{i}|a_{i},b_{i}\in L,m_{i}, n_{i}\in Z\}$.
Then by Eq. (2.3), one has
\begin{eqnarray*}
&&\sum_{i}(a_{i}\alpha^{n_{i}}\b^{m_i}([x,y]))b_{i}\\
&=&\sum_{i}[a_{i}\alpha^{n_{i}}\b^{m_i}([x,y]),b_{i}]\\
&+&\sum (R^2r^2 \widehat{R}^2\c \a^{-1}\b(b_{i}))(R^1\c \a\b^{-1}(a_{i})\alpha^{n_{i}}\b^{m_i}([r^1\c x,\widehat{R}^{1}\c y]))\\
&=&\sum_{i}[\alpha^{n_{i}}\b^{m_i}(\alpha^{-n_{i}}\b^{-m_i}(a_{i})[x,y]),b_{i}]\\
&+&\sum\alpha^{n_{i}+1}(((\alpha^{-n_{i}-1}(R^2r^2 \widehat{R}^2\c\alpha^{-1}\b(b_{i})\a\b^{-1}(R^1 \c a_{i}))))\alpha^{n_{i}+1}\b^{m_i-1}[r^1\c x,\widehat{R}^1\c y])\\
&=&\sum_{i}[\alpha^{n_{i}}\b^{m_i}(1_A(\alpha^{-n_{i}-1}\b^{-m_i-1}(a_{i})[x,y])1_{A}),b_{i}]\\
&+&\sum\alpha^{n_{i}+1}\b^{-m_i}((((\alpha^{-n_{i}-1}\b^{-m_i-1}(\alpha^{-1}(( R^2r^2 \widehat{R}^2\c b_{i}))(R^1\c a_{i_{0}}))))[r^1\c x,\widehat{R}^1\c y])1_{A})\\
&&\in U+U\subseteq U,
\end{eqnarray*}
since $\alpha(U)\subseteq U, \b(U)\subseteq U.$
Moreover, $I\neq 0$, for otherwise $[x,y]\in L$ will generate a nilpotent $H$-BiHom-ideal of $L$.
\hfill $\square$
\medskip

\noindent{\bf Theorem 3.3.}
Assume that $(L,\a,\b)$ is an prime  unital  generalized BiHom-associative algebra,
and $U$ be a $H$-BiHom-Lie ideal of $L$
such that $[U,U]\neq 0.$
Then there exists a $H$-BiHom-ideal $I$ of $L$
such that $0\neq[I,L]\subseteq U.$
\smallskip

{\it Proof }  Define $N_{L}(U)=\{x\in L|[x,L]\subseteq U\}.$
 Note that $U\subseteq N_{L}(U)$,  we have $[N_{L}(U),L]\subseteq U\subseteq N_{L}(H)$ and also $N_{L}(U)$ is a $H$-BiHom-Lie ideal.
It is easy to see that $[N_{L}(U),N_{L}(U)]\supseteq[U,U]\neq 0.$
Applying Lemma 3.2 to $N_{L}(U)$,
we may find a non-zero $H$-BiHom-ideal $I$ of $L$ such that $I\subseteq N_{L}(U)$, i.e., $[I,L]\subseteq U.$

Now we prove $[I,L]\neq 0$.
If not, choose $x\in I$ and $l,m\in L$,
then by Lemma 3.1(2),
$\a(x)[l,\a(m)]=[xl,\alpha\b(m)]-[x, R^2\c \b(m)]R^1\c\alpha(l),$ since $1_{A}(x\alpha(l))\in I$,
it is not hard to check that  $x[l,m]=0,$
$[L,L]\subseteq Ann_{L}(I).$
It is not hard to show that $Ann_{L}(I)$ is a $H$-BiHom-ideal.
In fact, let $x\in Ann_{L}(I)$ and $z\in I$.
This implies that $Ann_{L}(I)$ is a ideal.
$Ann_{L}(I)$ is clearly a $H$-BiHom-ideal.
It follows that $Ann_{L}(I)$ is a $H$-BiHom-ideal.
The simplicity of $L$ gives $[L,L]=0,$ a contradiction.
\hfill $\square$
\medskip

\noindent{\bf Corollary 3.4.}
Let $(L, \a,\b)$ be a simple unital generalized BiHom-associative algebra.
If $U$ is a $H$-BiHom-Lie ideal with $[U,U]\neq 0$, then $[L,L]\subseteq U.$
\medskip

As usual, we define a sequence of $H$-BiHom-ideals (the derived series) by
$$L^{0}=L,L^{(1)}=[L,L],L^{(2)}=[L^{(1)},L^{(1)}],\cdots,L^{(i)}=[L^{(i-1)},L^{(i-1)}].$$
$L$ is said to be solvable if $L^{(n)}=0$ for some $n.$
For $U\subseteq L$, $<U>$ denotes the $H$-BiHom-associative subalgebra of $L$ generated by $U.$
\medskip

\noindent{\bf Lemma 3.5.}
Let $U$ be a $H$-BiHom-ideal of $A$ and $<U>$   defined as above.
Then
 $<U>$ is a $H$-BiHom-Lie ideal of $A$,
 $[<U>,A]\subseteq U.$
\smallskip

{\it Proof}  Straightforward.
\hfill $\square$
\medskip

\noindent{\bf Lemma 3.6.}
Let $(L, \a,\b)$ be an $H$-simple  unital generalized BiHom-associative algebra. Then

(1) If $L^{(2)}\neq 0$, then $L=<L^{(1)}>.$

(2) If $L^{(3)}\neq 0$, then $L^{(2)}=L^{(1)}.$
\smallskip

{\it Proof}
(1) Let $S=\langle L^{(1)}\rangle.$
By Lemma 5.6(1), $S$ is a $H$-BiHom-Lie ideal of $L$.
By hypothesis, $[S,S]\supseteq L^{(2)}\neq 0$.
Obviously,
$\alpha(S)\subseteq S$.
Then by Lemma 3.2, $S$ contains a non-zero $H$-BiHom-ideal of $L$.
So $S=L$ since $L$ is $H$-simple.
\smallskip

(2) Let $V=L^{(2)}.$ The BiHom-Jacobi identity
implies that $V$ is a $H$-BiHom-Lie ideal of $L^{(3)}=[V,V]\neq 0$.
Also it is easy to see that $\alpha(V)\subseteq V,$   there is a $H$-BiHom-ideal $I$ of $L$. Thus $V\supseteq[L,L]$.
Clearly, $V\subseteq[L,L]$. This finishes the proof.\hfill $\square$
\medskip

\noindent{\bf Lemma 3.7.}
Let $V$ be an $H$-BiHom-Lie ideal of $[L,L]$ and
 $T(V)=\{l\in L|[l,L]\subseteq V\}.$
Then $T(V)$ is an $H$-BiHom-associative subalgebra of $L$.
\smallskip

{\it Proof } Straightforward.\hfill $\square$
\medskip

\noindent{\bf Lemma 3.8.}
With the notations as above, then
$[V,V]\subseteq T(V)$ and $[[T(V),T(V)],L]\subseteq T(V)$.
\smallskip

{\it Proof } Straightforward.  \hfill $\square$
\medskip

\noindent{\bf Theorem 3.9.}
Let $L$ be a simple  unital generalized BiHom-associative algebra
and $V\subseteq [L,L]$ be a $H$-BiHom-Lie ideal of $[L,L]$ such that $V\neq[L,L].$
Then $V$ is a solvable generalized BiHom-Lie subalgebra of $[L,L]$.
Moreover $[V,V]$ is nilpotent.

{\it Proof }   Firstly if  $L^{(3)}=0$,
 it is easy to see that $V^{(3)}=0$,
 and so $V$ is a solvable Hom-Lie subalgebra of $[L,L]$.
And if $L^{(3)}\neq0$,
suppose $[L,L]\subseteq T(V)$, then  $[L,[L,L]]\subseteq [L,T(V)]\subseteq V.$
But then $L^{(2)}=[[L,L],[L,L]]\subseteq [T(V),T(V)]\subseteq V.$
By Lemma 3.6(2), $L^{(2)}=[L,L]$, and so $[L,L]\subseteq V$,
which contradicts the hypothesis $V\neq[L,L].$
Thus we may assume that $[l,m]\bar{\in}T(V)$ for some $l,m\in L.$
By Lemma 3.7, we have
$[V+T(V),V+T(V)]\subseteq [V,V]+[V,T(V)]+[T(V),T(V)]\subseteq T(V),$
Hence $V+T(V)\neq L$.  So $L([[T(V),T(V)],[T(V),T(V)$ $]]L)\subseteq V+T(V)$.
This contradicts $L$ being $H$-simple, unless $[[T(V),T(V)],[T(V),T(V)]]=0$,
i.e.,$T(V)^{(2)}=0.$
Also, since $[V,V]\subseteq T(V)$, it follows that $V^{(3)}=0.$
Thus $V$ is a solvable $H$-BiHom-Lie subalgebra of $[L,L]$.
Moreover $[V,V]$ is nilpotent.

\begin{center}
 {\bf ACKNOWLEDGEMENT}
 \end{center}

   The paper is  supported by  the NSF of China (Nos. 11761017 and 11801304),    Guizhou Provincial  Science and Technology  Foundation (No. [2020]1Y005)
   and the Anhui Provincial Natural Science Foundation (Nos. 1908085MA03 and 1808085MA14).

\renewcommand{\refname}{REFERENCES}

\end{document}